\documentclass[leqno,12pt]{article}
\usepackage{amssymb,amsfonts}
\usepackage{amsmath,latexsym}

\usepackage{amstext,epic,eepic,epsf,pslatex}
\usepackage{graphicx}

\newcommand{\bepr}{{\em Proof} } 
\newcommand{\enpr}{\hfill \rule{.5em}{.5em}}

\newcommand{\R}{{\mathbb R}}

\newcommand{\Tr}{\hbox{Tr\,}}

\def\XXint#1#2#3{{\setbox0=\hbox{$#1{#2#3}{\int}$ }
\vcenter{\hbox{$#2#3$ }}\kern-.6\wd0}}

\newtheorem{thm}{Theorem}[section] 
\newtheorem{lemma}{Lemma}[section]

\newtheorem{cor}{Corollary}[section]

\begin{document}

\title{Estimating the number and the strength of collisions in molecular dynamics}

\author{Denis Serre \\ \'Ecole Normale Sup\'erieure de Lyon\thanks{U.M.P.A., UMR CNRS--ENSL \# 5669. 46 all\'ee d'Italie, 69364 Lyon cedex 07. France. {\tt denis.serre@ens-lyon.fr}}}

\date{}

\maketitle

\begin{abstract}
We consider the motion of a finite though large number of particles in the whole space $\R^n$. Particles move freely until they experience pairwise collisions. We use our recent theory of divergence-controlled positive symmetric tensors in order to establish two estimates regarding the set of collisions. The only information needed from the initial data is the total mass and the total energy.
\end{abstract}

\section{Models of molecular dynamics}

We consider a set of $N$ identical particles of mass $m$, moving in the whole space $\R^n$. The coordinate in the physical space $\R^{1+n}$ is denoted $x=(t,y)$ where $t$ is the time and $y\in\R^n$ the position. 

In practice $n=3$ and $N$ is of the size of the Avogadro number, but the analysis below is valid in every dimension and for any cardinality. We shall think of the particles as spheres of radius $a>0$, so that a collision between two particles occurs when their centers $y_\alpha$ approach to a distance $2a$~:
$$|y_\beta-y_\alpha|=2a.$$ 
The incoming and outgoing velocities ($v_{\alpha,\beta}$ and $v_{\alpha,\beta}'$ respectively) satisfy
\begin{equation}
\label{eq:incv}
(v_\beta-v_\alpha)\cdot(y_\beta-y_\alpha)<0,\qquad (v_\beta'-v_\alpha')\cdot(y_\beta-y_\alpha)\ge0.
\end{equation}
Our assumptions are as follows:
\begin{itemize}
\item Each particle $P_\alpha$ has a finite internal energy $\epsilon_\alpha\ge0$ and a velocity $v_\alpha\in\R^n$. These parameters remain constant between consecutive collisions involving $P_\alpha$.

\item When $P_\alpha$ and $P_\beta$ collide, at a given time $t$, their parameters (velocity, internal energy) experience a jump, which obey to the conservation of momentum and energy:
\begin{eqnarray*}
v_\alpha'+v_\beta' & = & v_\alpha+v_\beta,\\ 
\frac m2\left(|v_\alpha'|^2+|v_\beta'|^2\right)+\epsilon_\alpha'+\epsilon_\beta' & = & \frac m2\left(|v_\alpha|^2+|v_\beta|^2\right)+\epsilon_\alpha+\epsilon_\beta.
\end{eqnarray*}
\item The collisions are friction-less, meaning that the jump of velocity is orthogonal to the common tangent space to the particles:
\begin{equation}
\label{eq:normalimp}
v_\alpha'-v_\alpha=v_\beta-v_\beta'\parallel y_\beta-y_\alpha.
\end{equation}
\end{itemize}
The trajectory of a given particle is a polygonal chain. The conservation of energy implies a bound for the emerging velocities:
$$|v_\alpha'|^2+|v_\beta')|^2\le |v_\alpha|^2+\frac2m\,\epsilon_\alpha+|v_\beta|^2+\frac2m\,\epsilon_\beta.$$
The assumptions made above cover two important cases. On the one hand that of hard spheres, for which there is no exchange of internal energy (we might as well assume that there is no internal energy at all) and therefore
$$|v_\alpha'|^2+|v_\beta'|^2=|v_\alpha|^2+|v_\beta|^2.$$
On the other hand, we have the model of {\em sticky} particles, for which
$$v_\alpha'=v_\beta'=\frac12(v_\alpha+v_\beta).$$

\bigskip

Two important quantities emerge from this considerations, namely the total mass $M=Nm$ and the total energy
$$E=\sum_\alpha\left(\frac m2\,|v_\alpha(t)|^2+\epsilon_\alpha(t)\right),$$
which do not depend on the instant at which they are computed.  The third conserved quantity (total momentum)
$$Q=m\sum_\alpha v_\alpha(t)$$
will not be used below, for several reasons. On the first leg its nature is vectorial, which makes hard its use for an estimate. Besides, it just vanishes if we choose an appropriate inertial frame, a flaw which does not occur to the mass or the energy. Finally, it can be estimated by Cauchy--Schwarz inequality, $|Q|\le\sqrt{2ME}$\,, and therefore the knowledge of $M$ and $E$ will always give a sufficient information.

\section{Results}

Our aim is to estimate the number of collisions during the whole history, given only $M$ and $E$. We shall partly achieve this goal, in the sense that we establish estimates of the form
$$\sum_{\rm coll.}F(v_\alpha,v_\beta,v_\alpha',v_\beta')$$
in terms of $m,M$ and $E$, where $F$ is some explicit non-negative function, and we sum over all collisions. We confess that these estimates get poorer when either $v_\beta-v_\alpha$ is small (slow collision), or when the vectors $v_\alpha,v_\beta,v_\alpha',v_\beta'$ are approximately coplanar ({\em grazing} collisions). We believe that these limitations are inherent to the context of molecular dynamics, and that they do not reveal a weakness of the mathematical tool which we employ below.

Our main results are as follows. We limit ourselves to the case where every collision involves exactly two particles, and the evolution is defined for every time. This is a generic situation, as shown by Alexander in his PhD thesis \cite{Alex} ; see Theorem 4.2.1 of \cite{CIP}. See also Uchiyama's analysis \cite{Uch} of the Broadwell model, amodel with discrete velocity set. Our estimates involve the wedge product of two or three vectors in arbitrary dimensions~; we explain this notion, that comes from exterior calculus, in Paragraph \ref{ss:ext}. The product of $\ell$ vectors in $\R^n$ always vanishes if $\ell> n$.
\begin{thm}\label{th:main}
Consider a finite system of particles moving in the physical space $\R^n$ according to the laws described above. Let $M=Nm$ be the total mass and $E<\infty$ be the total energy of the system. Let us make the generic assumption that the collision set is locally finite, and that the motion involves only binary collisions.

At every collision, denote $v,v_1$ the incoming velocities and $v',v_1'$ the ougoing ones (they vary from one collision to another one, although the notation remains the same).

Then there exists a universal constant such that the following inequalities hold true
\begin{equation}
\label{eq:estbin}
m^2\sum_{\rm coll.}\left(\frac{E\,|v'-v|^2+M\,|v\wedge v'|^2}{\sqrt{(E+M\,|v|^2)(E+M\,|v'|^2)}}+\frac{E\,|v_1'-v_1|^2+M\,|v_1\wedge v_1'|^2}{\sqrt{(E+M\,|v_1|^2)(E+M\,|v_1'|^2)}}\right)\le c_nME,
\end{equation}
where the sum runs over the set of collisions. When $n\ge2$, we also have
\begin{equation}
\label{eq:esttri}
m^{\frac32}\sum_{\rm coll.}\frac{AB}C
\le c_nM^{\frac12}E^{\frac34}\,,
\end{equation}
where for each collision,
\begin{eqnarray*}
A & = & \left(E|(v'-v)\wedge(v_1-v)|^2+M|v\wedge v'\wedge v_1|^2\right)^{\frac12} \\
B & = & \left(4E+M(|v|^2+|v'|^2+|v_1|^2+|v_1'|^2)\right)^{1/4} \\
C & = &  \left((E+M|v|^2)(E+M|v'|^2)(E+M|v_1|^2)(E+M|v_1'|^2)\right)^{\frac14}.
\end{eqnarray*}
\end{thm}

By using the inequality $|v\wedge v'|\le|v|\cdot|v'-v|$, we see that (\ref{eq:estbin}) implies a simpler estimate:
\begin{cor}
\label{c:splr}
With the same assumptions as in Theorem \ref{th:main}, we have
\begin{equation}
\label{eq:estsplr}
m^2\sum_{\rm coll.}\frac{|v\wedge v'|^2}{|v|\cdot|v'|}\le c_nME
\end{equation}
for some universal constant $c_n$.
\end{cor}
It seems that (\ref{eq:esttri}) cannot be simplified in a similar way.

\paragraph{Comments.}
\begin{itemize}
\item These estimates are independent of the size $a$ of the particles. They are also dimension-independent, apart for the constants $c_n$.
\item The expressions $A,B$ and $C$ are symmetric functions of the velocities $v,v',v_1$ and $v_1'$ (obvious for $B$ and $C$).
\item The estimates above do not tell us whether the number of collisions is finite.
\item Nethertheless (\ref{eq:estbin}) tells us something about the number of the "strongest" ones. In terms of a typical velocity ${\bf v}:=\sqrt{E/M}\,$, the number of collisions for which $|v\wedge v'|^2\ge\epsilon{\bf v}^2|v|\cdot|v'|$ is bounded by $c_nN^2\epsilon^{-1}$.
\item Likewise, (\ref{eq:esttri}) tells us that the number of collisions for which $|(v'-v)\wedge(v_1-v)|\ge\epsilon{\bf v}^2$ and $|v|,|v'|,|v_1|,|v_1'|\le A{\bf v}$ for some $0<\epsilon<A<\infty$, is bounded by $C_{\epsilon,A}N^{\frac32}$.
On the contrary, it says very little about grazing collisions.
\end{itemize}

These polynomial bounds should be compared with those known for the whole set of collisions. For the hard spheres model with elastic collisions, Alexander \cite{Alex} and Vaserstein \cite{Vas} showed that this set is finite, a problem raised by Ya. Sinai. The proof of finiteness was simplified by Illner \cite{Illn,Illde}. So far, these works argued by contradiction and did not give an upper bound of the collision number. Later on Burago \& all. \cite{BFK} established the upper bound
$$(32N^{\frac32})^{N^2}.$$
This seems way too large to be accurate, since the fate of the particles is to disperse towards infinity in independent directions.

On the opposite side, some authors constructed configurations for which the number of collisions is superlinear in the number of particles. When the energy is not conserved and may increase arbitrarily, one can even observe infinitely many collisions \cite{Illn,ChIl}. For hard spheres with elastic collisions, Burdzy \& Duarte \cite{BD} anounce that some configuration leads to at least
$$\frac{N^3}{27}$$
collisions, a number which sounds reasonable. Let us point out that, according to our estimates, most of these collisions must be weak and grazing.

\paragraph{Plan of the paper.}
We begin by constructing in Section \ref{s:MMT} the mass-momentum tensor associated with the motion. It is a map $x\mapsto T(x)$ taking values in the cone of positive semi-definite matrices ${\bf Sym}_{1+n}^+$, albeit a singular one: the entries are bounded measures and their support is one-dimensional. The conservation of mass and momentum is expressed by the row-wise identity ${\rm Div}\,T=0$. One striking feature in this construction is the introduction of massless virtual particles ({\em collitons}) whose role is to carry the exchange of momentum between colliding particles. Section \ref{s:tools} presents the tools that will be used in the analysis. It is mainly a recall of our theory of Compensated Integrability for symmetric positive tensors whose row-wise divergence is a (vector-valued) bounded measure. We fix also a few notations related to exterior calculus. Sections \ref{s:bin} and \ref{s:ter} present the proofs of the binary and ternary estimates, respectively. We combine Compensated Integrability with a trick which we developed for the first time in \cite{Ser_X}.

\paragraph{Acknowledgement.} I am indebted to Laure Saint-Raymond and Reinhard Illner for valuable discussions and their help in gathering the relevant literature.

\section{The mass-momentum tensor}\label{s:MMT}

From now on, we denote $d=1+n$ the time-space dimension.

\subsection{Single particle}

We begin by considering a single particle $P$ whose constant velocity is $v\in\R^n$. The trajectory $t\mapsto (t,y(t))$ of the center of mass in the physical space $\R^{1+n}$ is a line $L$, whose direction is
$$\xi=\frac V{|V|}\,,\qquad\hbox{where } V:=\binom1v.$$
Let us define a symmetric matrix, whose entries are locally finite measures over $\R^{1+n}$, by
$$S=m\, V\otimes\xi\,\delta_L=m|V|\,\xi\otimes\xi\,\delta_L.$$
In other words
$$\langle S_{ab},\phi\rangle=mV_a\xi_b\int_\R\phi(\bar x+s\xi)\,ds,\qquad\forall\phi\in C_K(\R^{1+n}),$$
where $\bar x$ is any point on the line $L$.
\begin{lemma}
One has
$${\rm Div}\,S=0,$$
an equation that stands row-wise. More generally, if $Q\in \R^d$ with $Q\ne0$, and $\eta=\frac Q{|Q|}$\,, then for every line $L=\bar x+\R\eta$, the symmetric tensor
$$S^Q:=Q\otimes \eta\,\delta_L$$
 is divergence-free.
 \end{lemma}
 
 Notice that for a particle, $S$ is nothing but $S^Q$ with $Q=mV$. Remark also that $S^Q$ is everywhere positive semi-definite.
 
 \bigskip
 
 \bepr
 
 If $\phi$ is a test function, then
 $$\langle{\rm Div}\,S^Q,\phi\rangle=-\langle S^Q,\nabla\phi\rangle=-Q\int_\R \eta\cdot\nabla\phi(\bar x+s\eta)\,ds=-Q\int_\R \frac d{ds}\,\phi(\bar x+s\eta)\,ds=0.$$
 
 \enpr

\subsection{Multi-line configuration}

When $L$ is replaced by a semi-infinite line $L^+=\bar x+\R_+\eta$, the tensor
$$S^{Q+}:=Q\otimes \eta\,\delta_{L^+}$$
is no longer divergence-free. The calculation above yields
\begin{equation}
\label{eq:HL}
\langle{\rm Div}\,S^{Q+},\phi\rangle=\phi(\bar x)Q,
\end{equation}
which is recast as
$${\rm Div}\,S^{Q+}=Q\,\delta_{\bar x}\,.$$

Now, if finitely many vectors $Q_1,Q_2,\ldots$ are given, together with a point $\bar x\in\R^d$, we may form the converging semi-lines $L_j^+=\bar x +\R_+Q_j$ and define a symmetric tensor
$$S_{\rm mult}:=\sum_jS^{Q_j+}.$$
Then (\ref{eq:HL}) tells us that $S_{\rm mult}$ is divergence-free whenever
\begin{equation}
\label{eq:trpl}
\sum_jQ_j=0.
\end{equation}

\paragraph{Application to $1$-D molecular dynamics.} When $n=1$, we may simplify the model by setting $a=0$. 
At a binary collision, the point particles meet at some point $\bar x\in\R^{1+1}$, with incoming velocities $v,w$ and outgoing ones $v',w'$. Let us choose 
\begin{equation}
\label{eq:VVVV}
V_1=-\binom1v,\qquad V_2=-\binom1{w},\qquad V_3=\binom1{v'},\qquad V_4=\binom1{w'}
\end{equation}
and $Q_j=mV_j$. Then the positive semi-definite tensor
$$S^{Q_1+}+S^{Q_2+}+S^{Q_3+}+S^{Q_4+}$$
associated with this pair of particles is divergence-free~; the compatibility condition $Q_1+Q_2+Q_3+Q_4=0$ is ensured by the conservation of mass and momentum through the collision. Its support is the union of the trajectories.

\subsection{Binary collisions ($n\ge2$)}\label{ss:bin}

When $n\ge2$ instead, the radius $a$ must be positive, in order that collisions take place. 

Let two particles $P_i$ and $P_j$ collide at some time $t^*$. The trajectory of $P_i$ displays a kink at a point $\bar x_i=(t^*,\bar y_i)$, and that of $P_j$ does at $\bar x_j=(t^*,\bar y_j)$ at the same instant $t^*$. We have $|\bar y_j-\bar y_i|=2a$. Let us define $V_1,\ldots,V_4$ as in (\ref{eq:VVVV}). Locally, the trajectories are made of segments of  the semi-lines
$$L_1^+=\bar x_i+\R_+V_1,\qquad L_2^+=\bar x_j+\R_+V_2,\qquad L_3^+=\bar x_i+\R_+V_3,\qquad L_4^+=\bar x_j+\R_+V_4.$$
Let us denote again $Q_j=mV_j$. Because these do not meet at a single point, the tensor $S=S^{Q_1+}+S^{Q_2+}+S^{Q_3+}+S^{Q_4+}$ is not divergence-free. We have instead
$${\rm Div}\,S=(Q_1+Q_3)\delta_{x_i}+(Q_2+Q_4)\delta_{x_j}=(Q_1+Q_3)(\delta_{x_i}-\delta_{x_j}).$$

In order to recover a divergence-free tensor, we introduce a vector $Q$
$$Q=\binom0q,\qquad q=m(v'-v)=m(w-w').$$
Because of (\ref{eq:normalimp}), the segment $C=[\bar x_i,\bar x_j]$ has direction $Q$. In the neighbourhood of the collision, we can define the tensor 
$$T=S^{Q_1+}+S^{Q_2+}+S^{Q_3+}+S^{Q_4+}+Q\otimes\eta\,\delta_{C}.$$
Each of the five terms in the sum above is divergence-free away from either $\bar x_i$ or $\bar x_j$. At $\bar x_i$, ${\rm Div}\,T$ is a sum of three Dirac masses, whose weight is
\begin{equation}
\label{eq:VVQ}
m(V_1+V_3)-Q=\binom0{-mv+mv'-q}=0,
\end{equation}
where the minus sign in front of $Q$ comes from the fact that $Q$ is oriented from $x_j$ to $x_i$.
A similar identity holds true at $\bar x_j$, though with a plus sign. We conclude that 
$${\rm Div}\,T=0.$$

We may interpret the contribution $Q\otimes\eta\,\delta_C$ as that of a virtual particle. This particle is  mass-less, because the first component of $Q$ vanishes. It carries a momentum which is exchanged instantaneously between $P_i$ and $P_j$. We suggest the name {\em colliton} for this object. If  we took in account relativistic effects, there would be instead a pair of virtual particles (a particle and its anti-particle), travelling at the speed of light.

\subsection{The complete construction}

Assuming again that only binary collisions happen, we consider the union of trajectories of the centers of the $N$ particles. Each trajectory is a polygonal chain whose kinks occur where and when the particle suffers a collision. Each segment $J$ of a trajectory between two consecutive collisions contributes, as explained above, with the tensor
$$mV\otimes\xi\,\delta_J,\qquad V=\binom1v,\qquad\xi=\frac V{|V|}\,,$$
where $v$ is the particle velocity along $J$. At a collision we also have the contribution of the corresponding colliton, as described in the previous paragraph. The sum $T$ of all these contributions is a divergence-free positive semi-definite tensor, which we call the {\em mass-momentum} tensor of the configuration.

\bigskip

We point out that the support of $T$ is a graph, a one-dimensional object in $\R^{1+n}$. Thus $T$ vanishes almost everywhere in the Lebesgue sense. The support can be equiped with the positive measure $\Tr\,T$, with respect to which $T$ is rank-$1$ almost everywhere. We point out however that in order to apply Compensated Integrability (this tool is described in the next Section), we need to work with tensors of full rank $1+n$ (rather than rank-one), positive definite over a set of positive Lebesgue measure. To reach these goals, one option is to regularize $T$ by means of a convolution~; at a kink, this will allow us to combine three contributions, associated with $V_1,V_3$ and $Q$ with the notations of (\ref{eq:VVVV}) (two branches of a trajectory, plus the colliton). Such a combination has rank $2$ only since $m(V_1+V_3)=Q$. If we increase slightly the support of the convolution kernel, we may benefit of the combination of five contributions, those associated with $V_1,\ldots,V_4$ and $Q$, whose rank is generically $3$. This is still not sufficient if $n=3$ (the size of the mass-momentum tensors is $1+n$). We thus need an extra trick in order to conclude. The price to pay is that the resulting tensor will no longer be divergence-free~; this is the reason why we shall state Theorem \ref{th:CI} in the more general context of divergence-controlled tensors.

\section{Miscellaneous tools}\label{s:tools}

\subsection*{The determinant over ${\bf Sym}_d^+$}

The cone ${\bf Sym}_d^+$ of $d\times d$ positive semi-definite symmetric matrices is convex. On this cone, the map $K\mapsto\det K$ takes non-negative values, and it is monotonous: if $0_d\prec K\prec K'$ where the order is that between the associated quadratic forms, then
$$\det K\le\det K'.$$
Actually, the $d$-th root $K\mapsto(\det K)^{\frac1d}$ is a concave function over ${\bf Sym}_d^+$ (Theorem 6.10 in \cite{Ser_Mat}), homogeneous of degree one:
$$(K_1,\ldots,K_\ell\in{\bf Sym}_d^+)\Longrightarrow\left((\det(K+\cdots+K_\ell))^{\frac1d}\ge(\det K)^{\frac1d}+\cdots+(\det K_\ell)^{\frac1d}\right).$$ 
Since non-negative numbers $p_1,\ldots,p_\ell$ satisfy $(p_1+\cdots+p_\ell)^\alpha\ge p_1^\alpha+\cdots+p_\ell^\alpha$ for every exponent $\alpha\ge1$, we infer (take $\alpha=\frac d{d-1}$)
\begin{equation}
\label{eq:dmsun}
(K_1,\ldots,K_\ell\in{\bf Sym}_d^+)\Longrightarrow\left((\det(K+\cdots+K_\ell))^{\frac1{d-1}}\ge(\det K)^{\frac1{d-1}}+\cdots+(\det K_\ell)^{\frac1{d-1}}\right).
\end{equation}
Notice  that since this inequality is homogeneous of degree $\alpha>1$, we cannot deduce the concavity of $K\mapsto(\det K)^{\frac1{d-1}}$. The latter function is actually superlinear, thus non-concave.
\bigskip

Let $S$ be a symmetric tensor whose entries are distributions. Suppose that for every $\xi\in\R^d$, $\xi^TS\xi$ is non-negative, and therefore is a locally bounded measure. Then the entries of $S$ are locally bounded measures. In particular $\mu:=\Tr\,S$ is a non-negative locally bounded measure, and every $S_{ij}$ is absolutely continuous with respect to $\mu$. This allows us to define uniquely the expression $(\det S)^{\frac1d}$ as a locally bounded measure, absolutely continuous with respect to $\mu$~; just write $S_{ij}=s_{ij}\mu$ and define 
$$(\det S)^{\frac1d}=(\det s)^{\frac1d}\mu.$$ 
This definition inherits the properties of the determinant over ${\bf Sym}_d^+$, in particular the monotonicity and the concavity.

\subsection{Compensated integrability}

We shall make use of our recent theory of Compensated Integrability for divergence-controlled positive symmetric tensors, for which we refer to \cite{Ser_IHP,Ser_CI}. The appropriate version is given in the theorem below. We denote $\|\mu\|$ for the total mass of a (vector-valued) bounded measure $\mu$,
$$\|\mu\|=\langle|\mu|,{\bf1}\rangle.$$
This notation applies below in two distinct contexts, whether $\mu$ is a measure over an $(n+1)$-dimensional slab $H=(t_-,t_+)\times\R^n$, or a measure over $\R^n$. 
\begin{thm}\label{th:CI}
Let $H=(t_-,t_+)\times\R^n$ be a slab in $\R\times\R^n$, and $S$ be symmetric positive semi-definite tensor defined over $H$, whose entries are bounded measures. We assume that the row-wise divergence of $S$ is also a (vector-valued) bounded measure. Finally we assume that
the normal traces $S\vec e_t$ at the initial and final times $t=t_\pm$ are themselves bounded measures.

Then $(\det S)^{\frac1{n+1}}$ belongs to $L^{1+\frac1n}(H)$ and we have
\begin{equation}
\label{eq:CIn}
\int_H(\det S)^{\frac1n}dy\,dt\le c_n\left(\|S\vec e_t(t_-)\|+\|S\vec e_t(t_+)\|+\|{\rm Div}\,S\|\right)^{1+\frac1n},
\end{equation}
where $c_n$ is a finite constant independent of $S$ and $H$.
\end{thm}

\paragraph{Remarks.} 
\begin{itemize}
\item The assumption that ${\rm Div}\,S$ is a bounded measure allows us to define a normal trace in a rather weak space, here the dual of ${\rm Lip}(H)$. This is reminiscent to the definition of the normal trace of vector fields $\vec q\in H_{\rm div}(\Omega)$, used in the theory of incompressible fluids, which takes values in $H^{-1/2}(\partial\Omega)$.
\item The additional assumption that this trace is a bounded measure is equivalent to saying that  the extension $\tilde S$ by $0_{1+n}$ away from $H$ enjoys too the property that ${\rm Div}\,\tilde S$ is a bounded measure. Then
$${\rm Div}\,\tilde S=\widetilde{{\rm Div}\,S}-S\vec e_t(t_-)\otimes\delta_{t=t_-}+S\vec e_t(t_+)\otimes\delta_{t=t_+}$$
\item 
The qualitative part of the theorem is that the bounded measure $(\det S)^{\frac1d}$ is absolutely continuous with respect to the Lebesgue measure, and that its density is a function of class $L^{\frac d{d-1}}$, where $\frac d{d-1}=1+\frac1n$\,. The quantitative part (\ref{eq:CIn}) estimates this density.
\item
This theorem is useless when $S$ is rank-$1$ almost everywhere~; the estimate (\ref{eq:CIn}) is then trivial.
\end{itemize}

\subsection{Exterior calculus}\label{ss:ext}

We use a natural generalization of the determinant of a system of vectors. Let $Z^{(1)},\ldots,Z^{(k)}$ be a list of vectors in $\R^d$, with $k\le d$, and let ${\cal Z}\in{\bf M}_{d\times k}(\R)$ be the matrix whose columns are the $Z^{(j)}$'s. Let us decompose them on the canonical basis,
$$Z^{(j)}=\sum_{i=1}^dz_{ij}\vec e_i.$$
Then the exterior product $Z^{(1)}\wedge\cdots\wedge Z^{(k)}$ decomposes over the basis $\vec e_{i_1}\wedge\cdots\wedge\vec e_{i_k}$ with $i_1<\cdots<i_k$ of the exterior power $\Lambda^k(\R^d)$. 
Its coordinates are the $k\times k$ minors of $\cal Z$. For instance if $d=3$, $Z\wedge Z'$ identifies naturally with the cross product.

We denote
$$|Z^{(1)}\wedge\cdots\wedge Z^{(k)}|$$
the natural euclidian norm of the exterior product, which is invariant upon the action of the orthogonal group:
$$|RZ^{(1)}\wedge\cdots\wedge RZ^{(k)}|=|Z^{(1)}\wedge\cdots\wedge Z^{(k)}|,\qquad\forall R\in O_d(\R).$$
A pratical tool is given by the formula
$$|Z^{(1)}\wedge\cdots\wedge Z^{(k)}|^2=\det\left(\langle Z^{(j)}, Z^{(\ell)}\rangle\right)_{1\le j,\ell\le k}.$$
Suppose $k+\ell\le d$.
If ${\rm Span}(Z^{(1)},\ldots,Z^{(k)})$ and ${\rm Span}(W^{(1)},\ldots,W^{(\ell)})$ are orthogonal to each other, then
$$|Z^{(1)}\wedge\cdots\wedge Z^{(k)}\wedge W^{(1)}\wedge\cdots\wedge W^{(\ell)}|=|Z^{(1)}\wedge\cdots\wedge Z^{(k)}|\cdot|W^{(1)}\wedge\cdots\wedge W^{(\ell)}|.$$
In particular, we shall use the formula
$$|Z^{(1)}\wedge\cdots\wedge Z^{(k)}|=\left|\det(Z^{(1)},\ldots,Z^{(k)},v_{k+1},\ldots, v_d)\right|,$$
where $(v_{k+1},\ldots, v_d)$ is any unitary basis of the subspace ${\rm Span}(Z^{(1)},\ldots,Z^{(k)})^\bot$.

\section{The binary estimate}\label{s:bin}

In this section, we supplement the mass-momentum tensor $T$ with a contribution associated with the changes $v\mapsto v'$ in the direction of particles motions. Let $K$ be a kink of a trajectory, happening at a point $x^*$. The incoming/outgoing velocities being $v,v'$ respectively, with $v'\ne v$, we complete the free family
$$V=\binom1v,\quad V'=\binom1{v'}$$
into a basis $(V,V',z_2,\ldots,z_n)$ of $\R^{1+n}$, where $(z_2,\ldots,z_n)$ is some orthonormal basis of ${\rm Span}(V,V')^\bot$. To define the kink contribution, we introduce the positive semi-definite tensor
$$S_K=\sum_{j=2}^nz_j\otimes z_j\delta_{\sigma_j},$$
where $\sigma_j:=(x^*-az_j,x^*+az_j)$ is an interval.
Because ${\rm Div}\,S_K$ is a sum of Dirac masses at the end points $x^*\pm az_j$, we have 
$$\|{\rm Div}\,S_K\|=2(n-1).$$

Given a slab $H$, we form the tensor
$$T'=T+\sum_{\hbox{kinks in }H}b_KS_K,$$
where the positive numbers $b_K$ will be chosen later. We have
$$\|{\rm Div}\,T'\|\le2(n-1)\sum_{\hbox{kinks in }H}b_K,$$
where the inequality reflects the fact that a few end points could lie away from $H$.

Now, we make the convolution product $T_{2a}'=\phi_{2a}*T'$ with the non-negative kernel
$$\phi_{2a}=\frac1{|B_{2a}|}\,{\bf 1}_{B_{2a}}\,,$$
where $B_{2a}$ is the ball of radius $2a$ centered at the origin, and $\bf1$ denotes the characteristic function. This new tensor is still symmetric, positive semi-definite, and satisfies
$$\|{\rm Div}\,T_{2a}'\|\le2(n-1)\sum_{\hbox{kinks in }H}b_K.$$

Let $K$ be a kink in $H$, occuring at point $x^*$. If $x\in B_a(x^*)$, then $B_{2a}(x)$ contains $B_a(x^*)$. In particular $B_{2a}(x)$ meets every branch of the support of $T'$ around $x^*$ along a segment of length $\ge a$. There follows that
$$T_{2a}'(x)\ge\frac a{|B_{2a}|}\,\left(mV\otimes \xi+mV'\otimes \xi'+b_K\sum_2^nz_j\otimes z_j\right),\qquad\forall x\in B_a(x^*).$$
By the monotonicity of the determinant, we infer
$$\det T_{2a}(x)\ge\frac{m^2b_K^{n-1}}{a^{dn}|B_2|^d|V|\cdot|V'|}\,(\det(V,V',z_2,\ldots,z_n))^2=\frac{m^2b_K^{n-1}}{a^{dn}|B_2|^d}\,\frac{|V\wedge V'|^2}{|V|\cdot|V'|}\,.$$

We now apply Theorem \ref{th:CI} to $T_{2a}'$ in the slightly larger slab $H+B_{2a}$.
Let us assume first that the distance between two kinks in $H$ is $\ge2a$, so that the balls $B_a(x^*)$ are disjoint. Then the integral in the left hand side is estimated below by
\begin{equation}
\label{eq:below}
\int_{H+B_{2a}}(\det T_{2a}'(x))^{\frac1n}dx\ge\frac{2^{-d}m^{2/n}}{|B_2|^{1/n}}\,\sum_{\rm kinks}b_K^{1-\frac1n}\left(\frac{|V\wedge V'|^2}{|V|\cdot|V'|}\right)^{\frac1n}.
\end{equation}
When two (or more) balls $B_a(x^*)$ may overlap, we use Inequality (\ref{eq:dmsun}) to get
$$(\det T_{2a}(x))^{\frac1n}\ge\frac{m^{2/n}}{a^{d}|B_2|^{d/n}}\,\sum'b_K^{1-\frac1n}\left(\frac{|V\wedge V'|^2}{|V|\cdot|V'|}\right)^{\frac1n},$$
where the sum runs over all the kinks $K$ for which $x\in B_a(x^*)$. By integrating over $x\in H$ and then rearranging the sum, we deduce again the same lower bound (\ref{eq:below}).

From  (\ref{eq:CIn}) and (\ref{eq:below}), we deduce
$$m^{2/n}\sum_{\rm kinks}b_K^{1-\frac1n}\left(\frac{|V\wedge V'|^2}{|V|\cdot|V'|}\right)^{\frac1n}\le
c_n\left(\|T_{2a}'\vec e_t(t_--2a)\|+\|T_{2a}'\vec e_t(t_++2a)\|+2(n-1)\sum_{\rm kinks}b_K\right)^{1+\frac1n}$$
for some universal constant $c_n$~; we recall that the kinks in the summation are those in $H$. There remains to estimate the masses in the right hand side. We notice that
$$\|T_{2a}'\vec e_t(t_++2a)\|\le\sup_{t>t_+}\|T'\vec e_t(t)\|=\sup_{t>t_+}\|T\vec e_t(t)\|=m\sup_{t>t_+}\sum|V|,$$
where the sum runs over all particles. Because of $|V|=\sqrt{1+|v|^2}\le1+\frac12\,|v|^2$, we infer
$$\|T_{2a}'\vec e_t(t_++2a)\|\le M+E$$
and we deduce an explicit inequality
\begin{equation}
\label{eq:expli}
m^{2/n}\sum_{\rm kinks}b_K^{1-\frac1n}\left(\frac{|V\wedge V'|^2}{|V|\cdot|V'|}\right)^{\frac1n}\le
c_n\left(M+E+(n-1)\sum_{\rm kinks}b_K\right)^{1+\frac1n}
\end{equation}
with a slightly different constant $c_n$.

We exploit (\ref{eq:expli}) by letting the coefficients $b_K$ vary. Let us choose positive numbers $\beta_K$ and another positive parameter $\lambda$. By setting
$$b_K=\lambda\beta_K^{\frac n{n-1}},$$
we have
$$m^{2/n}\sum_{\rm kinks}\beta_K\left(\frac{|V\wedge V'|^2}{|V|\cdot|V'|}\right)^{\frac1n}\le
c_n\lambda^{\frac1n-1}\left(M+E+\lambda (n-1)\sum_{\rm kinks}\beta_K^{\frac n{n-1}}\right)^{1+\frac1n}.$$
Choosing
$$\lambda=\frac{M+E}{\sum\beta_K^{\frac n{n-1}}}\,,$$
we obtain
$$m^{2/n}\sum_{\rm kinks}\beta_K\left(\frac{|V\wedge V'|^2}{|V|\cdot|V'|}\right)^{\frac1n}\le
c_n(M+E)^{\frac2n}\left(\sum_{\rm kinks}\beta_K^{\frac n{n-1}}\right)^{1-\frac1n}=c_n(M+E)^{\frac2n}\left\|\vec\beta\right\|_{\ell^{n/(n-1)}}.$$
Since $\frac n{n-1}$ is the conjugate exponent of $n$, a convenient choice of the vector $\vec\beta$ gives us an estimate of the $\ell^n$-norm of the vector whose coordinates are the expressions
$$\left(\frac{|V\wedge V'|^2}{|V|\cdot|V'|}\right)^{\frac1n}.$$
Specifically, we have proved
\begin{equation}
\label{eq:elln}
m^2\sum_{\rm kinks}\frac{|V\wedge V'|^2}{|V|\cdot|V'|}\le c_n(M+E)^2,
\end{equation}
with again a slightly different constant $c_n$. This estimate can be recast in terms of the velocities only, by using
$$|V|=\sqrt{1+|v|^2}\,,\qquad |V\wedge V'|^2=|v-v'|^2+|v\times v'|^2.$$
We point out that the constant in (\ref{eq:elln}) depends only on the space dimension and not on the data. In particular, it does not depend upon the width of the slab $H$. Since the right-hand side depends only on the initial data, we deduce that (\ref{eq:elln}) is valid when we sum over the whole set of kinks in $\R^{1+n}$.

\bigskip

The last step involves a rescaling. Keeping the same space variable $y$ and changing the time variable into $\bar t=\mu^{-1}t$ has the effect of replacing the velocities $v$ by $\mu v$. The mass is preserved, while the energy becomes $\mu^2 E$. Applying (\ref{eq:elln}) to the new dependent/independent variables, we obtain a parametrized estimate
$$m^2\sum_{\rm kinks}\frac{\mu^2|v'-v|^2+\mu^4|v\wedge v'|^2}{\sqrt{(1+\mu^2|v|^2)(1+\mu^2|v'|^2)}}\le c_n(M+\mu^2E)^2,\qquad\forall\mu>0.$$
Choosing now $\mu^2=M/E$\,, we arrive at our final estimate, which is now invariant upon the rescaling of the time and space variables:
$$m^2\sum_{\rm kinks}\frac{E\,|v'-v|^2+M\,|v\wedge v'|^2}{\sqrt{(E+M\,|v|^2)(E+M\,|v'|^2)}}\le c_nME.$$

\section{The ternary estimate}\label{s:ter}

We change slightly the strategy described in the previous section, by considering the set of collisions instead of that of kinks. At a collision $C$ between two particles $P_\alpha$ and $P_\beta$, we denote $\bar x$ the middle point between the kinks $x_{\alpha,\beta}$. The incoming/outgoing velocities are denoted $v,v_1,v',v_1'$. When the vectors $v'-v$ and $y_\beta-y_\alpha$ are not colinear\footnote{This is the generic situation where the trajectories are not coplanar}, the vectors $V,V'$ and $x_\beta-x_\alpha$ are not coplanar~; they span a $3$-dimensional subspace which contains also $V_1$ and $V_1'$. Then we choose an orthonormal basis $(w_3,\ldots,w_n)$ of the subspace ${\rm Span}(V,V',x_\beta-x_\alpha)^\bot$, and we define a tensor 
$$S_C=\sum_{j=3}^nw_j\otimes w_j\,\delta_{\theta_j}\,,\qquad\theta_j:=(\bar x-aw_j,\bar x+aw_j).$$
As above, ${\rm Div}\,S_C$ is a sum of Dirac masses at the end points $\bar x\pm aw_j$, whose total mass is $2(n-2)$.

Let $H=(t_-,t_+)\times\R^n$ be a slab. We say that a collision $C$ is in $H$ if it happens at a time $t^*\in(t_-,t_+)$. We form the symmetric, positive semi-definite tensor
$$T''=T+\sum_{\hbox{coll. in }H}b_CS_C,$$
where the positive numbers $b_C$ are to be chosen. Then we make the convolution $T_{3a}''$ of $T''$ with the kernel $\phi_{3a}$. The total mass of ${\rm Div}\,T''$ is $2(n-2)\sum b_C$.

Consider a collision $C$  in $H$. For any point $x$ in $B_a(\bar x)$, the ball $B_{3a}(x)$ contains the balls $B_a(x_{\alpha,\beta})$, thus meets the five segments involved in the collision pattern along intervals of lengths at least $a$. Therefore we have
$$T_{3a}''(x)\ge\frac a{|B_{3a}|}\,\left(m(V\otimes \xi+V'\otimes \xi'+V_1\otimes \xi_1+V_1'\otimes\xi_1'+\bar V\otimes\bar\xi)+b_C\sum_3^nw_j\otimes w_j\right),\qquad\forall x\in B_a(\bar x)$$
where
$$\bar V=\binom0{v'-v},\qquad\xi=\frac{\bar V}{|\bar V|}$$
accounts for the colliton. With the monotonicity of the determinant, we infer the coarse lower bound
$$\det T_{3a}''(x)\ge\frac{m^3b_C^{n-2}}{a^{dn}|B_3|^d}\,\frac{|V\wedge V'\wedge V_1|^2}{|V|\cdot|V'|\cdot|V_1|}\,.$$
More generally, if $x$ is $a$-close to several collisions, then we have 
$$(\det T_{3a}''(x))^{\frac1n}\ge\frac{m^{3/n}}{a^{d}|B_3|^{d/n}}\,\sum'b_C^{1-2/n}\left(\frac{|V\wedge V'\wedge V_1|^2}{|V|\cdot|V'|\cdot|V_1|}\right)^{\frac1n}\,,$$
where the sum runs over those collisions $C$ in $H$ for which $x\in B_a(\bar x)$.

Integrating with respect to $x$ in $H+B_a$, and rearranging the sum, we obtain
$$\int_{H+B_a}(\det T_{3a}''(x))^{\frac1n}\,dx\ge\frac{m^{3/n}|B_1|}{|B_3|^{d/n}}\,\sum_{\hbox{coll. in }H}b_C^{1-2/n}\left(\frac{|V\wedge V'\wedge V_1|^2}{|V|\cdot|V'|\cdot|V_1|}\right)^{\frac1n}\,,$$
Applying (\ref{eq:CIn}) to $T_{3a}''$ in $H+B_a$, we deduce again
$$m^{\frac3n}\sum_{\hbox{coll. in }H}b_C^{1-2/n}\left(\frac{|V\wedge V'\wedge V_1|^2}{|V|\cdot|V'|\cdot|V_1|}\right)^{\frac1n}\le c_n\left(M+E+(n-2)\sum_{\hbox{coll. in }H}b_C\right)^{1+\frac1n}.$$
Introducing as above auxiliary positive parameters $\lambda$ and $\beta_C$, with
$$b_C=\lambda\beta_C^p,\qquad p:=\frac n{n-2}\,,$$
and choosing
$$\lambda=\frac{M+E}{\beta_C^p}\,,$$
we deduce
$$m^{\frac3n}\sum_{\hbox{coll. in }H}\beta_C\left(\frac{|V\wedge V'\wedge V_1|^2}{|V|\cdot|V'|\cdot|V_1|}\right)^{\frac1n}\le c_n(M+E)\|\vec\beta\|_{\ell^p}.$$
An appropriate choice of $\vec\beta$ yields an estimate in the $\ell^{p'}$-norm, where $p'=\frac n2$\,. Namely, we have
\begin{equation}
\label{eq:nonsym}
m^{\frac32}\sum_{\hbox{coll. in }H}\frac{|V\wedge V'\wedge V_1|}{\sqrt{|V|\cdot|V'|\cdot|V_1|}}
\le c_n(M+E)^{\frac32}.
\end{equation}

Remark that we used only the three first factors among $V\otimes \xi,V'\otimes  \xi',V_1\otimes  \xi_1$ and $V_1'\otimes  \xi_1'$ in order to bound by below the determinant of $T_{3a}''(x)$. We might as well have taken any three of the four. Remarking that the wedge product remains the same, up to the sign, because of $V+V_1=V'+V_1'$, we deduce that (\ref{eq:nonsym}) remains valid when we replace $|V|\cdot|V'|\cdot|V_1|$ by any other product of three norms in the denominator. In other terms, we have
$$m^{\frac32}\sum_{\hbox{coll. in }H}\frac{|V\wedge V'\wedge V_1|}{\sqrt{|V|\cdot|V'|\cdot|V_1|\cdot|V_1'|}}\,\left(|V|^{1/2}+|V'|^{1/2}+|V_1|^{1/2}+|V_1'|^{1/2}\right)
\le c_n(M+E)^{\frac32},$$
where we now sum a quantity that is symmetric in the four velocities.

Our next remark is, as before, that the left-hand side does not depend upon the choice of the slab. Thus our inequality is 
valid when we sum over the set of all collisions in $\R^{1+n}$. 

We end our analysis as above by applying the rescaling argument. This transform an inhomogeneous estimate into an homogeneous one:
$$m^{\frac32}\sum_{\rm coll.}\frac{\left(E|(v'-v)\wedge(v_1-v)|^2+M|v\wedge v'\wedge v_1|^2\right)^{\frac12}\left(4E+M(|v|^2+|v'|^2+|v_1|^2+|v_1'|^2)\right)^{1/4}}{\left((E+M|v|^2)(E+M|v'|^2)(E+M|v_1|^2)(E+M|v_1'|^2)\right)^{\frac14}}\,
\le c_nM^{\frac12}E^{\frac34}\,,$$
where we have used also the identity
$$|V\wedge V'\wedge V_1|^2=|(v'-v)\wedge(v_1-v)|^2+|v\wedge v'\wedge v_1|^2.$$

\end{document}